\documentclass[12pt,a4paper,reqno]{amsart}
\usepackage{amssymb}
\usepackage{dcpic,pictex}
\usepackage{multirow}
\usepackage{graphicx}
\usepackage{slashbox}
\usepackage{cite}

\usepackage{hyperref}
\usepackage[main=english,french]{babel}
\newtheorem{definition}{Definition}[section]

\newtheorem{theorem}{Theorem}[section]
\newtheorem{remark}{Remark}[section]
\newcommand{\D}{\mathbb{D}}

\newcommand{\I}{\mathbb{I}}
\newcommand{\R}{\mathbb{R}}

\newcommand{\Com}{\mathbb{C}}
\newcommand{\N}{\mathbb{N}}

\begin{document}


\title[General fractional integrals and derivatives]{General fractional integrals and 
        derivatives with the Sonine kernels}

\author{Yuri Luchko}
\curraddr{Beuth Technical University of Applied Sciences Berlin,  
     Department of  Mathematics, Physics, and Chemistry,  
     Luxemburger Str. 10,  
     13353 Berlin,
Germany}
\email{luchko@beuth-hochschule.de}

\subjclass[2010]{26A33; 26B30; 44A10; 45E10}
\dedicatory{}
\keywords{Sonine kernel, general fractional derivative, general fractional integral, $n$-fold general fractional derivative, $n$-fold general fractional integral, fundamental theorem of Fractional Calculus}

\begin{abstract}
In this paper, we address the general fractional integrals and derivatives with the Sonine kernels on the spaces of functions with an integrable singularity at the point zero. First, the Sonine kernels and their important special classes and particular cases are discussed. In particular, we introduce a class of the Sonine kernels that possess an integrable singularity of power function type at the point zero. For the general fractional integrals and derivatives with the Sonine kernels from this class, two fundamental theorems of  fractional  calculus are proved. Then,  we construct the $n$-fold general fractional integrals and derivatives and study their properties. 
\end{abstract}

\maketitle


\section{Introduction}\label{sec:1}

Starting from the paper \cite{Koch11}, a considerable interest in Fractional Calculus (FC) research  has been  devoted to the general fractional derivatives (GFD) of convolution type (see also \cite{Han20,LucYam20}) defined in the Riemann--Liouville sense
\begin{equation}
\label{FDR-L}
(\D_{(k)}\, f) (t) = \frac{d}{dt}\int_0^t k(t-\tau)f(\tau)\, d\tau
\end{equation}
and in the 
Caputo sense
\begin{equation}
\label{FDC}
( _*\D_{(k)}\, f) (t) = (\D_{(k)}\, f) (t) - f(0)k(t).
\end{equation}

The Riemann--Liouville and the Caputo fractional derivatives of order $\alpha,\ 0< \alpha < 1$ are particular cases of the GFDs \eqref{FDR-L} and \eqref{FDC} with the kernel
\begin{equation}
\label{single}
k(t) = h_{1-\alpha}(t),\ 0 <\alpha <1,\ h_{\beta}(t) := \frac{t^{\beta -1}}{\Gamma(\beta)},\ \beta >0.
\end{equation}
They are defined as follows:
\begin{equation}
\label{RLD}
(D^\alpha\, f)(t) = \frac{d}{dt}(I^{1-\alpha}_{0+} f)(t),\ t>0,
\end{equation}
\begin{equation}
\label{CD}
( _*D^\alpha\, f)(t) = \frac{d}{dt}(I^{1-\alpha}_{0+} f)(t) - f(0)h_{1-\alpha}(t),\ t>0, 
\end{equation}
\vspace{6pt}
$I^\alpha_{0+}$ being the Riemann--Liouville fractional integral of order $\alpha\ (\alpha >0)$:
\begin{equation}
\label{RLI}
\left(I^\alpha_{0+} f \right) (t)=\frac{1}{\Gamma (\alpha )}\int\limits_0^t(t -\tau)^{\alpha -1}f(\tau)\,d\tau,\ t>0.
\end{equation}

Other particular cases of the GFDs \eqref{FDR-L} and \eqref{FDC} are the multi-term fractional derivatives and the fractional derivatives of the
distributed order. They are generated by \eqref{FDR-L} and \eqref{FDC} with the kernels 
\begin{equation}
\label{multi}
k(t) = \sum_{k=1}^n a_k\, h_{1-\alpha_k}(t),
\ \ 0 < \alpha_1 <\dots < \alpha_n < 1,\ a_k \in \R,\ k=1,\dots,n, 
\end{equation}
\begin{equation}
\label{distr}
k(t) = \int_0^1 h_{1-\alpha}(t)\, d\rho(\alpha),
\end{equation}
respectively, where $\rho$ is a Borel measure defined on the interval $[0,\, 1]$.

One of the main issues addressed in \cite{Koch11} is to clarify the conditions which allow an interpretation of the integro-differential operators \eqref{FDC} as a kind of the fractional derivatives. In particular, these derivatives and the appropriate defined fractional integrals have to satisfy the fundamental theorem of FC (see \cite{Luc20} for a discussion of this theorem). Moreover, the solutions to the time-fractional differential equations of certain types with the GFDs are expected to behave as the ones to the evolution equations. This includes the property of complete monotonicity of solutions to the relaxation equation and the positivity of the fundamental solution to the Cauchy problem for the fractional diffusion equation with the time-derivatives of this type. In \cite{Koch11}, a class of kernels of the GFD \eqref{FDC} that satisfy the requirements formulated above  is described in terms of their Laplace transforms. 

In the meantime, some other results regarding the ordinary and partial time-fractional differential equations with the GFDs of Caputo type have been already derived. The ordinary fractional differential equations with the Caputo type GFDs   are  addressed in \cite{Koch11,Koch19_2,KK,Sin18} and the time-fractional diffusion equations with these derivatives are  considered in \cite{Koch11,Koch19_1,Koch19_2,LucYam16,LucYam20,Sin20}. The publications in \cite{JK17,KJ19_1,KJ19} are devoted to the inverse problems for the GFDs and for the partial differential equations containing the general fractional integrals and derivatives. 

In this paper, we address the general fractional integrals and derivatives with the Sonine kernels that possess an integrable singularity of power function type at the point zero. In particular, the $n$-fold general fractional integrals and derivatives are constructed and studied for the first time. Even if we restrict ourselves to the mathematical notions, results, and proofs, the constructions presented in this paper are potentially useful for applications. They provide a wider scope for fitting the kernels of the fractional derivatives to the data available from the measurements. The optimal fitting of the GFDs kernels to the data at hand can be realised by solving the inverse problems for the models in form of the fractional differential equations with the GFDs (see, e.g.,     \cite{JK17,KJ19_1,KJ19} for some results of this kind). Moreover, nowadays,  the multi-term fractional differential equations with the Caputo, Riemann--Liouville, or Hilfer fractional derivatives are actively employed for modelling of complex systems and processes with memory that possess different time scales. The multi-term fractional differential equations with the $n$-fold GFDs is a far reaching generalisation of these models. These equations and their potential applications will be considered elsewhere. 

The rest of the paper is organised as follows. In  Section \ref{sec:2}, the basic spaces of functions that are employed for derivation of our main results are shortly discussed. Section \ref{sec:3} addresses the Sonine condition and the Sonine kernels. We provide both some examples of the Sonine kernels and descriptions of several important classes of the Sonine kernels. In particular, we introduce a class of the Sonine kernels that possess an integrable singularity of power function type at the point zero. In Section 4, the general fractional integrals (GFI) and derivatives with the Sonine kernels from this class are introduced and studied on the spaces of functions that can have an integrable singularity at the point zero. The first and   second fundamental theorems of FC for the GFDs are proved on the appropriate spaces of functions. Moreover, the $n$-fold general fractional integrals and derivatives that correspond to the Riemann--Liouville and Caputo derivatives of an arbitrary order are constructed and their basic properties are studied.


\section{Preliminaries}
\label{sec:2}

In this section, we introduce the spaces of functions that are used in the further discussions and present their main properties. 
\begin{definition}
\label{dd1}
By $C_{-1}(0,+\infty)$, we denote the space of functions continuous on the positive real semi-axis that can have an integrable singularity of a power function type at the point zero: 
\begin{equation}
\label{2-7}
C_{-1}(0,+\infty)\, = \, \{f:\ f(t)=t^{p-1}f_1(t),\ t>0,\ p>0,\ f_1\in C[0,+\infty)\}.
\end{equation}
\end{definition}
A family of the spaces $C_\alpha(0,+\infty),\ \alpha \ge -1$ was first introduced by  \cite{Dim66} for construction of an operational calculus for the hyper-Bessel differential operator and then employed for the operational calculi for the fractional derivatives in \cite{BasLuc95, Luc99, LucSri95, LucYak94, YakLuc94} and other related publications. 

Using the well-known properties of the Laplace convolution
\begin{equation}
\label{2-2}
(f\, *\, g)(t) = \int_0^{t}\, f(t-\tau)g(\tau)\, d\tau,
\end{equation}
the following theorem can be easily proved:

\begin{theorem}[\cite{LucGor99}]
\label{t1}
The triple $\mathcal{R}_{-1} = (C_{-1}(0,+\infty),+,*)$ with the usual
addition $+$ and   multiplication $*$ in the form of the Laplace convolution \eqref{2-2} is a commutative ring without divisors of zero.
\end{theorem}

In particular, let us mention that the integration operator 
\begin{equation}
\label{2-3}
(I^1_{0+}\, f)(t) = (\{1\}\, *\, f)(t) = \int_0^t\, f(\tau)\, d\tau
\end{equation}
can be interpreted as a multiplication in the ring $\mathcal{R}_{-1}$. In the representation \eqref{2-3}, by $\{1\}$,  we denote  the function that is identically equal to $1$ on $[0,+\infty)$. 

Because of the inclusion $C[0,+\infty) \subset C_{-1}(0,+\infty)$ and because a sum and a convolution of any two continuous functions are again continuous functions, we have the inclusion
\begin{equation}
\label{2-8}
(C[0,+\infty),+,*) \subset (C_{-1}(0,+\infty),+,*). 
\end{equation}

In the next sections, the sub-spaces $C_{-1}^m(0,+\infty),\ m\in \N_0=\N \cup\{0\}$ of the space $C_{-1}(0,+\infty)$   are  also used. They are defined as follows: A function $f=f(t),\ t>0$, is said to belong to the space $C_{-1}^m(0,+\infty)$ if and only if $f^{(m)} \in C_{-1}(0,+\infty)$. 

The spaces $C_{-1}^m(0,+\infty)$ were first introduced and studied by \cite{LucGor99}. Here, we  only  present some relevant properties of these spaces  (for the proofs,  see \cite{LucGor99}). 

\begin{theorem}[\cite{LucGor99}]
\label{t2}
For the family of the spaces $C_{-1}^m(0,+\infty),\ m\in \N_0$, the following statements hold true:
\vskip 0.1cm 

\noindent
(1) $C_{-1}^0(0,+\infty)\, \equiv \, C_{-1}(0,+\infty)$.
\vskip 0.1cm 

\noindent
(2) $C_{-1}^m(0,+\infty), \ m\in \N_0$ is a vector space over the field $\R$ (or $\Com$).
\vskip 0.1cm 

\noindent
(3) If $f\in C_{-1}^m(0,+\infty)$ with $m\ge 1$, 
then $f^{(k)}(0+) := \lim\limits_{t\to 0+} f^{(k)}(t) <+\infty,\ 0\le k\le 
m-1$, and the function
$$
\tilde f (t) = 
\begin{cases}f(t), & t>0, \\
f(0+), & t=0
\end{cases}
$$
belongs to the space $C^{m-1}[0,+\infty)$.

\vskip 0.1cm 

\noindent
(4) If $f\in C_{-1}^m(0,+\infty)$ with $m\ge 1$, then $f \in 
C^{m}(0,+\infty)\cap C^{m-1}[0,+\infty)$.
\vskip 0.1cm 

\noindent
(5) For $ m\ge 1$, the following representation holds true:
$$
f\in C_{-1}^m(0,+\infty) \Leftrightarrow
f(t) = (I^m_{0+} \phi)(x) + \sum \limits_{k=0}^{m-1} f^{(k)}(0)\frac{t^k}{k!},\ t\ge 
0,\ \phi \in C_{-1}(0,+\infty).
$$
\vskip 0.1cm 

\noindent
(6) Let $f\in C_{-1}^m(0,+\infty), \ m\in \N_0$,
$f(0)=\dots=f^{(m-1)}(0)=0$ and $g\in C_{-1}^1$. Then,  the 
Laplace convolution $h(t) = (f*g)(t)$
belongs to the space $C_{-1}^{m+1}(0,+\infty)$ and $h(0)=\dots=h^{(m)}(0)=0$.
\end{theorem}

Theorem \ref{t2} ensures that the triples $\mathcal{R}_{-1}^m = (C_{-1}^m(0,+\infty),+,*)$ are sub-rings of the basic ring $\mathcal{R}_{-1}$. 

In the rest of this section, we list some examples of the functions from the space $C_{-1}(0,+\infty)$ and the formulas for their multiplication in the ring $\mathcal{R}_{-1}$ (Laplace convolution).

We start with the important power functions (see \eqref{single} for their definition):
\begin{equation}
\label{2-9}
(h_{\alpha} \, * \, h_\beta)(t) \, = \, h_{\alpha+\beta}(t),\ \alpha,\beta >0,\ t>0.
\end{equation}

For $\alpha,\beta >0$, the functions $h_{\alpha},\ h_{\beta}$, and $h_{\alpha+\beta}$ all belong to the space $C_{-1}(0,+\infty)$.   Formula \eqref{2-9} is an easy consequence from the well-known representation of the Euler Beta-function in terms of the Gamma-function:
$$
B(\alpha,\beta) := \int_0^1 (1-\tau)^{\alpha -1}\, \tau^{\beta -1}\, d\tau \, = \, \frac{\Gamma(\alpha)\Gamma(\beta)}{\Gamma(\alpha+\beta)},\ \alpha,\beta>0.
$$ 

The same formula can be applied to verify a more general relation 
\begin{equation}
\label{2-10}
(h_{\alpha,\rho} \, * \, h_{\beta,\rho})(t) \, = \, h_{\alpha+\beta,\rho}(t),\ \alpha,\beta >0,\ \rho \in \R,\ t>0,
\end{equation}
where the function $h_{\alpha,\rho} \in C_{-1}(0,+\infty)$ is given by the expression
\begin{equation}
\label{2-11}
h_{\alpha,\rho}(t) = \frac{t^{\alpha -1}}{\Gamma(\alpha)}\, e^{-\rho t},\ \ \alpha >0,\ \rho \in \R,\ t>0.
\end{equation}

Because $\mathcal{R}_{-1}$ is a ring, we can define the integer order convolution powers of a function $\kappa \in C_{-1}(0,+\infty)$ as follows:
\begin{equation}
\label{2-12}
\kappa^n :=\underbrace{\kappa*\ldots\ * \kappa}_n,\ n\in \N.
\end{equation}

For example,   by      repeatedly applying   Formula \eqref{2-9}, we get the representation
\begin{equation}
\label{2-13}
h_{\alpha}^n(t) = h_{n\alpha}(t),\ n\in \N.
\end{equation}

In particular, we mention the well-known formula
\begin{equation}
\label{2-13-1}
\{1\}^n(t) = h_{1}^n(t)= h_{n}(t)=\frac{t^{n-1}}{\Gamma(n)} = \frac{t^{n-1}}{(n-1)!},\ n\in \N.
\end{equation}




\section{The Sonine Kernels}
\label{sec:3}

In the short notices \cite{Son} published in 1884, Sonine addressed a generalisation of the Abel integral equation and applied Abel's method for its analytical solution. He recognised that the most essential ingredient of Abel's solution is the formula
\begin{equation}
\label{3-1}
(h_{\alpha}\, * \, h_{1-\alpha})(t) = h_1(t) = \{1 \},\ 0<\alpha<1,\ t>0
\end{equation}
that is a particular case of \eqref{2-9}. Note that, for $0<\alpha<1$, the inclusions $h_{\alpha} \in C_{-1}(0,+\infty)$ and $h_{1-\alpha} \in C_{-1}(0,+\infty)$ are valid. 

As a generalization of \eqref{3-1}, Sonine suggested to consider the functions $\kappa$ and $k$ that satisfy the relation
\begin{equation}
\label{3-2}
(\kappa \, *\, k )(t) = \{1 \},\ t>0.
\end{equation}

Nowadays, such functions are called the Sonine kernels and the relation \eqref{3-2} is called the Sonine condition. For a Sonine kernel $\kappa$, the kernel $k$ that satisfies the Sonine condition \eqref{3-2} is called a dual or associate kernel to $\kappa$. Of course, $\kappa$ is then an associated kernel to $k$. In what follows, we denote the set of the Sonine kernels by $\mathcal{S}$.

Evidently, not all  functions  defined on $\R_+$  are  a Sonine kernel. For example,         the constant function $\{ 1\}$ does not possess any associate kernel. The problem of analytical description of the set $\mathcal{S}$ of the Sonine kernels is still not completely solved. However, several important classes of the Sonine kernels have been introduced in the literature and in this section some relevant results are shortly presented. 

In \cite{Son}, Sonine addressed the case of the kernels that can be represented in the form
\begin{equation}
\label{3-3}
\kappa(t) = h_{\alpha}(t) \cdot \, \kappa_1(t),\ \kappa_1(t)=\sum_{k=0}^{+\infty}\, a_k t^k, \ a_0 \not = 0,\ 0<\alpha <1,
\end{equation}
where the function $\kappa_1=\kappa_1(t)$ is analytical on $\R$. He showed that such functions are always the Sonine kernels and their (unique) associate kernels can be represented as follows:
\begin{equation}
\label{3-4}
k(t) = h_{1-\alpha}(t) \cdot k_1(t),\ k_1(t)=\sum_{k=0}^{+\infty}\, b_k t^k. 
\end{equation}

The coefficients $b_k,\ k\in \N_0$ are uniquely determined by the coefficients $a_k,\ k\in \N_0$ as solutions to the following triangular system of linear equations:
\begin{equation}
\label{3-5}
a_0b_0 = 1,\ \sum_{k=0}^n\Gamma(k+1-\alpha)\Gamma(\alpha+n-k)a_{n-k}b_k = 0,\ n\ge 1.
\end{equation}

In \cite{Son}, some examples of the kernels from $\mathcal{S}$ in the form \eqref{3-3} and \eqref{3-4} were derived including the prominent pair
\begin{equation}
\label{Bess}
\kappa(t) = (\sqrt{t})^{\alpha-1}J_{\alpha-1}(2\sqrt{t}),\ 
k(t) = (\sqrt{t})^{-\alpha}I_{-\alpha}(2\sqrt{t}),\ 0<\alpha <1,
\end{equation}
where 
$$
J_\nu (t) = \sum_{k=0}^{+\infty} \frac{(-1)^k(t/2)^{2k+\nu}}{k!\Gamma(k+\nu+1)},\ 
I_\nu (t) = \sum_{k=0}^{+\infty} \frac{(t/2)^{2k+\nu}}{k!\Gamma(k+\nu+1)}
$$
are the Bessel and   modified Bessel functions, respectively. 

Another example of this type is the following pair of the Sonine kernels 
(\cite{Zac08}):
\begin{equation}
\label{3-6} 
\kappa(t) = 
h_{\alpha,\rho}(t),\ \rho\ge 0,\ 0<\alpha <1,
\end{equation}
\begin{equation}
\label{3-7} 
k(t) = h_{1-\alpha,\rho}(t) \, +\, \rho\, \int_0^t h_{1-\alpha,\rho}(\tau)\, d\tau,\ t>0,
\end{equation}
where the function $h_{\alpha,\rho}$ is defined by the formula \eqref{2-11}. 

In \cite{Wick}, the condition of analyticity of the function $\kappa_1$ from \eqref{3-3}  is replaced by a weaker condition that the series for $\kappa_1$ has a finite convergence radius $r>0$. Then,  the function $\kappa$ defined by \eqref{3-3} is still a Sonine kernel and its unique associate kernel is given by \eqref{3-4}, where the series for $k_1$ has the same convergence 
radius $r$. 

Of course, the class $\mathcal{S}$ of the Sonine kernels is not restricted to the functions in form \eqref{3-3} and includes, for instance, some functions with the power-logarithmic singularities at the origin (\cite{Rub, Sam}). 

A completely different approach for description of the Sonine kernels was suggested by \cite{Koch11}, who established  a connection between the complete Bernstein functions and the Sonine kernels. More precisely, they showed  that any function $k = k(t),\ t>0$ that satisfies Conditions (K1)--(K4) listed below is a Sonine kernel.

\begin{enumerate}

\item[(K1)] The Laplace transform $\tilde k$ of $k$,
\begin{equation}
\label{Laplace} 
\tilde k(p) = ({\mathcal L}\, k)(p)\ =\ \int_0^{+\infty} k(t)\, e^{-pt}\, dt
\end{equation}
exists for all $p>0$.

\item[(K2)] $\tilde k(p)$ is a Stieltjes function (see,  e.g., \cite{[SSV]} for definition and properties of the Stieltjes functions).

\item[(K3)] $\tilde k(p) \to 0$ and $p \tilde k(p) \to +\infty$ as $p \to +\infty$.

\item[(K4)] $\tilde k(p) \to +\infty$ and $p \tilde k(p) \to 0$ as $p \to 0$.

\end{enumerate}


In what follows, we denote the set of the kernels that satisfy Conditions (K1)--(K4) by $\mathcal{K}$. As   mentioned above,  $\mathcal{K}\subset \mathcal{S}$. Moreover, the solutions to the fractional differential equations with the GFDs of type \eqref{FDC} with the kernels $k \in \mathcal{K}$ behave as the ones to the evolution equations  \cite{Koch11}. In particular, the unique solution to the Cauchy problem for the fractional relaxation equation with a positive initial condition is complete monotone and the fundamental solution to the Cauchy problem for the fractional diffusion equation with the time-derivative of this type can be interpreted as a probability density function. 

The Sonine kernels from $\mathcal{K}$  are  introduced in terms of their Laplace transforms and using the Sonine condition \eqref{3-2} in the Laplace domain. Providing the Laplace transforms $\tilde \kappa,\ \tilde k$ of the functions $\kappa$ and $k$   exist, the convolution theorem for the Laplace transform leads to the equation
\begin{equation}
\label{Laplace-Sonine} 
\tilde \kappa(p) \cdot \tilde k(p) = \frac{1}{p},\ \Re(p)>p_{\kappa,k} \in \R
\end{equation}
for the Laplace transforms of the Sonine kernels $\kappa$ and $k$. However, the Laplace transforms of the Sonine kernels do not always exist. One can easily construct such kernels based on the Sonine Formula \eqref{3-3} with the analytical functions $\kappa_1$ that grow  more quickly  than any exponential function $\exp(ct),\, c >0$ as $t\to +\infty$. One of the examples of this sort is the Sonine kernel
 \begin{equation}
\label{3-3-Lap}
\kappa(t) = h_{\alpha}(t) \cdot \exp(t^2),\ 0<\alpha <1.
\end{equation}

In \cite{Han20}, an important class of the Sonine kernels  is introduced in terms of the completely monotone functions. Namely, it is shown there that any singular (unbounded in a neighborhood of the point zero) locally integrable
completely monotone function is a Sonine kernel. We denote this class of the Sonine kernels by $\mathcal{H}$ ($\mathcal{H} \subset \mathcal{S}$). Moreover, if $\kappa \in \mathcal{H}$, then its associate kernel $k$ also does belong to $\mathcal{H}$. A typical example of a pair of the Sonine kernels from $\mathcal{H}$ is as follows (\cite{Han20}):
\begin{equation}
\label{3-8} 
\kappa(t) = h_{1-\beta+\alpha}(t)\, +\, h_{1-\beta}(t),\ 0<\alpha < \beta <1,
\end{equation}
\begin{equation}
\label{3-9} 
k(t) = t^{\beta -1}\, E_{\alpha,\beta}(-t^\alpha),
\end{equation}
where $E_{\alpha,\beta}$ stands for the two-parameters Mittag--Leffler function that is defined by the following convergent series:
\begin{equation}
\label{ML}
E_{\alpha,\beta}(z) = \sum_{k=0}^{+\infty} \frac{z^k}{\Gamma(\alpha\, k + \beta)},\ \alpha >0,\ \beta,z\in \Com.
\end{equation}

To prove this, the well-known Laplace transform formula 
$$ 
\left({\mathcal L}\, t^{\beta-1}E_{\alpha,\beta}(-t^\alpha)\right)(p)\, = \, \frac{p^{\alpha-\beta}}{p^\alpha +1},\ \alpha,\, \beta >0
$$
for the two-parameters Mittag--Leffler function is applied along with   Formula 
\eqref{Laplace-Sonine} and the known fact that the function \eqref{3-8} is a singular locally integrable
completely monotone function (\cite{[SSV]}).

As mentioned in \cite{Sam} (see also \cite{Han20}), any Sonine kernel  possesses an integrable singularity at the point zero. On the other hand, the kernels of the fractional integrals and derivatives should be singular (\cite {DGGS}). Thus, the fractional integrals and derivatives with the Sonine kernels are worth   being  investigated. 

In the next section, we deal with the general fractional integrals and derivatives with the Sonine kernels that belong to the space $C_{-1}(0,+\infty)$. Because the Sonine kernels have an integrable singularity at the point zero, we restrict ourselves to the corresponding sub-space of $C_{-1}(0,+\infty)$:
\begin{equation}
\label{subspace}
 C_{-1,0}(0,+\infty) \, = \, \{f:\ f(t) = t^{p-1}f_1(t),\ t>0,\ 0<p<1,\ f_1\in C[0,+\infty)\}.
\end{equation}
\begin{definition}
\label{dd2}
Let $\kappa,\, k \in C_{-1,0}(0,+\infty)$ be a pair of the Sonine kernels, i.e., the Sonine condition \eqref{3-2} be fulfilled. The set of such Sonine kernels    is denoted by $\mathcal{S}_{-1}$:
\begin{equation}
\label{Son}
(\kappa,\, k \in \mathcal{S}_{-1} ) \ \Leftrightarrow \ (\kappa,\, k \in C_{-1,0}(0,+\infty))\wedge ((\kappa\, *\, k)(t) \, = \, \{1\}).
\end{equation}
\end{definition}
Because the ring $\mathcal{R}_{-1}$ is divisors free (Theorem \ref{t1}), for any Sonine kernel from $\mathcal{S}_{-1}$, its associate kernel is unique. It is also worth mentioning that the kernels of the conventional time-fractional derivatives and integrals introduced so far do belong to the set $\mathcal{S}_{-1}$ and thus the constructions presented in the next section are applicable for both   the known and   many new FC operators. In particular, all examples of the Sonine kernels provided in this section including the general kernels \eqref{3-3} and \eqref{3-4} introduced by Sonine in his original publication \cite{Son} are from the set $\mathcal{S}_{-1}$. 


\section{Fractional Integrals and Derivatives with the Sonine Kernels}
\label{sec:4}



In this section, we address the GFIs and the GFDs with the Sonine kernels from $\mathcal{S}_{-1}$ on the spaces of functions $C_{-1}^m(0,+\infty), \ m\in \N_0$ discussed in Section \ref{sec:2} and on their sub-spaces. 
\begin{definition}
\label{dd3}
Let $\kappa \in \mathcal{S}_{-1}$. The GFI with the kernel $\kappa$ is defined by the formula
\begin{equation}
\label{GFI}
(\I_{(\kappa)}\, f)(t) := (\kappa\, *\, f)(t) = \int_0^t \kappa(t-\tau)f(\tau)\, d\tau,\ t>0.
\end{equation}
\end{definition}
The GFI \eqref{GFI} with the kernel $\kappa(t) = h_\alpha(t) =\frac{t^{\alpha-1}}{\Gamma(\alpha)}\in \mathcal{S}_{-1},\ 0<\alpha<1$ is the Riemann--Liouville fractional integral \eqref{RLI}:
\begin{equation}
\label{GFI_1}
(\I_{(\kappa)}\, f)(t) = (I^\alpha_{0+}\, f)(t),\ t>0.
\end{equation} 

For the Riemann--Liouville fractional integral, the relations 
\begin{equation}
\label{GFI_11}
(I^0_{0+}\, f)(t) = f(t),\ \ (I^1_{0+}\, f)(t) = \int_0^t\, f(\tau)\, d\tau 
\end{equation} 
hold true. Combining the Formulas \eqref{GFI_1} and \eqref{GFI_11}, it is natural to define the GFIs with the kernels $\kappa(t) = h_0(t)$ and $\kappa(t) = h_1(t)$ (that do not belong to the set $\mathcal{S}_{-1}$) as follows:
\begin{equation}
\label{GFI_12}
(\I_{(h_0)}\, f)(t) := (I^0_{0+}\, f)(t) = f(t),\ \ (\I_{(h_1)}\, f)(t) := (I^1_{0+}\, f)(t) = \int_0^t\, f(\tau)\, d\tau. 
\end{equation} 

In the definition \eqref{GFI_12}, the function $h_0$ is interpreted as a kind of the $\delta$-function that plays the role of a unity with respect to multiplication in form of the Laplace convolution. In particular, one could extend   Formula \eqref{3-1} valid for $0<\alpha<1$ in the usual sense to the case $\alpha =1$ that has to be understand in the sense of generalised functions:
\begin{equation}
\label{3-1-1}
(h_{1}\, * \, h_{0})(t) := h_1(t) = \{1 \},\ t>0.
\end{equation}

Using the known Sonine kernels, many other particular cases of the GFI \eqref{GFI} can be constructed. Here,  we mention just a few of them:
\begin{equation}
\label{GFI_2}
(\I_{(\kappa)}\, f)(t) = (I^{1-\beta+\alpha}_{0+}\, f)(t) + (I^{1-\beta}_{0+}\, f)(t) ,\ t>0
\end{equation} 
with the Sonine kernel $\kappa$ defined by \eqref{3-8},
\begin{equation}
\label{GFI_3}
(\I_{(\kappa)}\, f)(t) = \frac{1}{\Gamma(\alpha)}\int_0^t (t-\tau)^{\alpha -1}\, e^{-\rho (t-\tau)}\, f(\tau)\, d\tau,\ t>0
\end{equation} 
with the Sonine kernel $\kappa$ defined by \eqref{3-6}, and 
\begin{equation}
\label{GFI_4}
(\I_{(\kappa)}\, f)(t) = \int_0^t (\sqrt{t-\tau})^{\alpha-1}J_{\alpha-1}(2\sqrt{t-\tau})f(\tau)\, d\tau,\ t>0
\end{equation}
with the Sonine kernel $\kappa$ defined by \eqref{Bess}.

Regarding the properties of the GFI \eqref{GFI} on the space $C_{-1}(0,+\infty)$, they can be easily derived from the known properties of the Laplace convolution and Theorem \ref{t2}. In particular, in what follows, we need the mapping property
\begin{equation}
\label{GFI-map}
\I_{(\kappa)}:\, C_{-1}(0,+\infty)\, \rightarrow C_{-1}(0,+\infty),
\end{equation}
the commutativity law 
\begin{equation}
\label{GFI-com}
\I_{(\kappa_1)}\, \I_{(\kappa_2)} = \I_{(\kappa_2)}\, \I_{(\kappa_1)},\ \kappa_1,\, \kappa_2 \in \mathcal{S}_{-1},
\end{equation}
and the index law
\begin{equation}
\label{GFI-index}
\I_{(\kappa_1)}\, \I_{(\kappa_2)} = \I_{(\kappa_1*\kappa_2)},\ \kappa_1,\, \kappa_2 \in \mathcal{S}_{-1}
\end{equation}
that are valid on the space $C_{-1}(0,+\infty)$. 

\begin{definition}
\label{dd4}
Let $\kappa \in \mathcal{S}_{-1}$ and $k$ be its associate Sonine kernel. 

The GFDs of the Riemann--Liouville and the Caputo types associate to the GFI \eqref{GFI} are defined as follows (see \eqref{FDR-L} and \eqref{FDC}): 
\begin{equation}
\label{GFDL}
(\D_{(k)}\, f) (t) = \frac{d}{dt}(k\, * \, f)(t),\ t>0,
\end{equation}
\begin{equation}
\label{GFDC}
( _*\D_{(k)}\, f) (t) = (\D_{(k)}\, f) (t) - f(0)k(t),\ t>0.
\end{equation}
\end{definition}

It is easy to see that the GFD \eqref{GFDC} in the Caputo sense can be represented as a regularized GFD \eqref{GFDL} in the Riemann--Liouville sense:
\begin{equation}
\label{GFDC_new}
( _*\D_{(k)}\, f) (t) = (\D_{(k)}\, [f(\cdot)-f(0)]) (t),\ t>0.
\end{equation}

In this paper, we mainly deal with the GFDs in form \eqref{GFDC} on the spaces $C_{-1}^m(0,+\infty),$ $m\in \N$ and their sub-spaces. For the functions from $C_{-1}^1(0,\infty)$, Theorem \ref{t2} and the known formula for differentiation of the integrals depending on parameters lead to an important representation of the GDF \eqref{GFDL}:
\begin{equation}
\label{GFDL-1}
(\D_{(k)}\, f) (t) = (k\, * \, f^\prime)(t) + f(0)k(t),\ t>0.
\end{equation}

Thus, for $f\in C_{-1}^1(0,+\infty)$, the GFD \eqref{GFDC} can be rewritten in the form
\vspace{6pt}
\begin{equation}
\label{GFDC_1}
( _*\D_{(k)}\, f) (t) = (k\, * \, f^\prime)(t),\ t>0. 
\end{equation}

It was known already to Sonine that the operators of type \eqref{GFDL} are left inverse to the GFI \eqref{GFI}. However, he provided just a formal derivation that needs a justification for the given sets of the kernels and the spaces of functions. In our case, we do this in the following theorem (see \cite{Koch11} for the case of the GFDs with the kernels $k\in \mathcal{K}$ and \cite{Han20} for the case $k\in \mathcal{H}$).

\begin{theorem}[First Fundamental Theorem for the GFD]
\label{t3}
Let $\kappa \in \mathcal{S}_{-1}$ and $k$ be its associate Sonine kernel. 

Then,  the GFD \eqref{GFDL} is a left inverse operator to the GFI \eqref{GFI} on the space $C_{-1}(0,+\infty)$: 
\begin{equation}
\label{FTL}
(\D_{(k)}\, \I_{(\kappa)}\, f) (t) = f(t),\ f\in C_{-1}(0,+\infty),\ t>0,
\end{equation}
and the GFD \eqref{GFDC} is a left inverse operator to the GFI \eqref{GFI} on the space $C_{-1,k}(0,+\infty)$: 
\begin{equation}
\label{FTC}
( _*\D_{(k)}\, \I_{(\kappa)}\, f) (t) = f(t),\ f\in C_{-1,k}(0,+\infty),\ t>0,
\end{equation}
where $C_{-1,k}(0,+\infty) := \{f:\ f(t)=(\I_{(k)}\, \phi)(t),\ \phi\in C_{-1}(0,+\infty)\}$.
\end{theorem}
\begin{proof}
The proof of    Formula \eqref{FTL} follows the lines of the Sonine derivations:
$$
(\D_{(k)} I_{(\kappa)}\, f) (t) = \frac{d}{dt}(k * \kappa * f)(t) = 
\frac{d}{dt}(\{ 1\} * f)(t) = \frac{d}{dt} \int_0^t f(\tau)\, d\tau = f(t),\ t>0.
$$

Because $k,\, \kappa \in \mathcal{S}_{-1}$ and $f\in C_{-1}(0,+\infty)$, the  validity of the formulas in the chain of the equations above is justified by Theorem \ref{t2}. 


To prove    Formula \eqref{FTC}, we first show that any function $f$ from the space $C_{-1,k}(0,+\infty)$ satisfies the properties
\begin{equation}
\label{Ck}
\I_{(\kappa)} f \in C_{-1}^1(0,+\infty)\ \mbox{ and }\ (\I_{(\kappa)}\, f)(0) = 0.
\end{equation}

Indeed, the following chain of implications holds true ($t>0$):
$$
f=\I_{(k)}\, \phi = k*\phi \ \Rightarrow \I_{(\kappa)} f = \kappa * f = \kappa * k *\phi = \{ 1 \} *\phi \ \Rightarrow 
$$
$$
(\I_{(\kappa)} f)^\prime = (\{ 1 \} *\phi)^\prime = \phi \in C_{-1}(0,+\infty).
$$

To prove that $(\I_{(\kappa)} f)(0) = 0$, we use the representation from the previous formula and estimate the function $\I_{(\kappa)} f$ on the interval $(0,1]$ ($\phi(t) = t^{p-1}\phi_1(t),\ p>0,\ \phi_1\in C[0,+\infty)$):
$$
|(\I_{(\kappa)} f)(t)| = |(\{ 1 \} *\phi)(t)| \le \int_0^t |\phi(t)|dt = \int_0^t \tau^{p-1}|\phi_1(\tau)|d\tau \le 
$$
$$
C \int_0^t \tau^{p-1}d\tau = 
C \frac{t^p}{p},\ p>0,\ t\in (0,\, 1],\ C\in \R_+.
$$

Because the function $\I_{(\kappa)} f \in C_{-1}^1(0,+\infty)$ is continuous on $[0,+\infty)$ (Theorem \ref{t2}), this estimate immediately leads to the desired result:
$$
(\I_{(\kappa)} f)(0) \, = \, \lim_{t\to 0+0} (\I_{(\kappa)} f)(t) = 0.
$$

Now,  we employ    Formula \eqref{FTL} that is already proved above (evidently, $C_{-1,k}(0,+\infty) \subset C_{-1}(0,+\infty)$) and finalise the proof:
$$
( _*\D_{(k)}\, \I_{(\kappa)}\, f) (t) = (\D_{(k)}\, \I_{(\kappa)}\, f) (t) - (\I_{(\kappa)}\, f) (0)\, k(t) = f(t),\ t>0.
$$
\end{proof}

\begin{remark}
\label{r1}
In the proof of Theorem \ref{t3}, we show the implication
$$
(f = \I_{(k)}\phi, \phi \in C_{-1}(0,+\infty)) \Rightarrow (\I_{(\kappa)} f \in C_{-1}^1(0,+\infty)\, \wedge \, (\I_{(\kappa)}\, f)(0) = 0).
$$

The inverse implication holds true, too:
$$
(\I_{(\kappa)} f \in C_{-1}^1(0,+\infty)\, \wedge \, (\I_{(\kappa)}\, f)(0) = 0) \Rightarrow (f = \I_{(k)}\phi, \phi \in C_{-1}(0,+\infty)).
$$

Indeed, $\I_{(\kappa)} f \in C_{-1}^1(0,+\infty)$ implicates the inclusion $(\I_{(\kappa)} f)^\prime = (\kappa * f)^\prime = \phi \in C_{-1}(0,+\infty)$.
Integrating this formula, we arrive at the representation
$$
\{1\}* (\kappa * f)^\prime = \kappa * f = \{1\} * \phi 
$$
because $(\I_{(\kappa)}\, f)(0) = (\kappa * f)(0) = 0$. Then,  we multiply the last relation by $k$:
$$
k*\kappa * f = \{ 1\} * f = k* \{1\} * \phi = \{1\} *(k* \phi), \ t>0.
$$

Differentiation of the relation $\{ 1\} * f = \{1\} *(k* \phi)$ leads to the representation \mbox{$f = \I_{(k)}\phi$}, $\phi \in C_{-1}(0,+\infty)$.

Thus,  we prove  that 
$$
\{f:\ f(t)=(\I_{(k)}\, \phi)(t),\ \phi\in C_{-1}(0,+\infty)\} \, = 
$$
$$
 \{f:\ \I_{(\kappa)} f \in C_{-1}^1(0,+\infty)\, \wedge \, (\I_{(\kappa)}\, f)(0) = 0\}.
$$

It is worth mentioning that this characterization of the space $C_{-1,k}(0,+\infty)$ is a generalization of the known property of the spaces of functions employed while treating the Abel integral equation (Theorem 2.3 in \cite{Samko}) 
$$ 
\left\{ f:\, f = I^{1-\alpha}\, \phi,\ \phi \in L_1(0,1)\right\}\, = \, 
\left\{ f:\, I^{\alpha}f\in \mbox{AC}([0,\, 1])\, \wedge \, (I^{\alpha}f)(0) = 0\right\}
$$
to the case of the GFI with the Sonine kernel $\kappa\in \mathcal{S}_{-1}$.
\end{remark}

Now,   we provide the explicit formulas for the GFDs $ _*\D_{(k)}$ that correspond to the particular cases \eqref{GFI_1} and  \eqref{GFI_2}--\eqref{GFI_4} of the GFI \eqref{GFI}. To make the formulas more compact, we employ the representation \eqref{GFDC_1}. However, these formulas with the evident modifications are also valid for the GFDs in form \eqref{GFDC}. 

The Sonine kernel associate to the kernel $\kappa(t) = h_\alpha(t),\ 0<\alpha<1$ is $k(t) = h_{1-\alpha}(t)$, and thus the GFD in form \eqref{GFDC_1} associate to the Riemann--Liouville integral \eqref{RLI} is the Caputo derivative
\begin{equation}
\label{GFD_1}
( _*\D_{(k)}\, f)(t)=( _*D^\alpha\, f)(t) = (I^{1-\alpha}_{0+}\, f^\prime)(t).
\end{equation} 

In    Formula \eqref{GFI_12}, the Riemann--Liouville fractional integral  is  defined for $\alpha=0$ and $\alpha = 1$ that corresponds to the kernels $\kappa(t) = h_0(t)$ 
and $\kappa(t) = h_1(t)$. The convolution Formula \eqref{3-1-1} suggests the following natural definition of the GFD with the kernels $k = h_1(t)$ 
and $k(t) = h_0(t)$, respectively:
\begin{equation}
\label{GFD_12}
( _*\D_{(h_1)}\, f)(t):=( _*D^0\, f)(t) = f(t),\ ( _*\D_{(h_0)}\, f)(t):=( _*D^1\, f)(t) = \frac{df}{dt}. 
\end{equation} 

Of course, to justify the definitions \eqref{GFD_12}, the suitable subsets of the GFD domain should be introduced, where the formulas
\begin{equation}
\label{GFD_13}
\lim_{k\to h_1}\, _*\D_{(k)}\, f = f,\ \lim_{k\to h_0}\, _*\D_{(k)}\, f = \frac{df}{dt} 
\end{equation} 
are valid in some sense (see \cite{HL19} for a description of these and other properties that the fractional derivatives and integrals should possess). This will be done elsewhere. 

The GFD that corresponds to the GFI \eqref{GFI_2} has the Mittag--Leffler function \eqref{3-9} in the kernel:
\begin{equation}
\label{GFD_2}
( _*\D_{(k)}\, f)(t) = \int_0^t (t-\tau)^{\beta -1}\, E_{\alpha,\beta}(-(t-\tau)^\alpha)\, f^\prime(\tau)\, d\tau,\ 0<\alpha < \beta <1,\ t>0.
\end{equation}

The Sonine kernel associate to the kernel $\kappa(t) = h_{\alpha,\rho}(t),\ 0<\alpha<1,\ \rho\ge 0$ is given by    Formula \eqref{3-7} and 
thus the GFD associate to the GFI \eqref{GFI_3} has the form 
\begin{equation}
\label{GFD_3}
( _*\D_{(k)}\, f)(t) = \int_0^t \left( h_{1-\alpha,\rho}(t-\tau) \, +\, \rho\, \int_0^{t-\tau} h_{1-\alpha,\rho}(\tau_1)\, d\tau_1 \right)\, f^\prime(\tau)\, d\tau,\ t>0.
\end{equation} 
Finally, the GFD that corresponds to the GFI \eqref{GFI_4} takes the form (see the Sonine pair \eqref{Bess}): 
\begin{equation}
\label{GFD_4}
( _*\D_{(k)}\, f)(t) = \int_0^t 
(\sqrt{(t-\tau)})^{-\alpha}I_{-\alpha}(2\sqrt{(t-\tau)})f^\prime(\tau)\, d\tau,\ t>0.
\end{equation}

Now,  we proceed with the 2nd fundamental theorem of FC for the GFDs in the Riemann--Liouville and Caputo senses (see \cite{Koch11} for the case of the GFDs with the kernels $k\in \mathcal{K}$ and \cite{Han20} for the case $k\in \mathcal{H}$).

\begin{theorem}[Second  Fundamental Theorem for the GFD]
\label{t4}
Let $\kappa \in \mathcal{S}_{-1}$ and $k$ be its associate Sonine kernel.

Then,  the relations
\begin{equation}
\label{2FTC}
(\I_{(\kappa)}\, _*\D_{(k)}\, f) (t) = f(t)-f(0),\ t>0,
\end{equation}
\begin{equation}
\label{2FTL}
(\I_{(\kappa)}\, \D_{(k)}\, f) (t) = f(t),\ t>0
\end{equation}
hold valid for the functions $f\in C_{-1}^1(0,+\infty)$.
\end{theorem}

\begin{proof}
As      shown  at the beginning of this section, for the functions 
$f\in C_{-1}^1(0,+\infty)$, the GFD $ _*\D_{(k)}$ can be represented in form \eqref{GFDC_1}. Moreover, the functions from $C_{-1}^1(0,+\infty)$ are continuous on $[0,\, +\infty)$ (Theorem \ref{t2}) and $f^\prime \in C_{-1}(0,+\infty)$. Then,  we have the following chain of equations:
$$
(\I_{(\kappa)}\, _*\D_{(k)}\, f) (t) = (\kappa * ( _*\D_{(k)}\, f))(t) = (\kappa * k*f^\prime)(t) = 
$$
$$
(\{ 1\} *f^\prime)(t) = \int_0^t f^\prime(\tau)d\tau = f(t)-f(0).
$$

For the GFD in the Riemann--Liouville sense, we first employ    Formula \eqref{GFDL-1} and then proceed as follows:
$$
(\I_{(\kappa)}\, \D_{(k)}\, f) (t) = (\I_{(\kappa)}\, ((k*f^\prime)(\tau) + f(0)k(\tau))) (t) =
$$
$$
(\I_{(\kappa)}\, _*\D_{(k)}\, f) (t) + f(0)(\kappa * k)(t) = f(t)-f(0) + f(0) = f(t).
$$
\end{proof}

In the rest of this section, we consider the compositions of the GFIs and construct the suitable fractional derivatives in the sense of the 1st Fundamental Theorem of FC. 

\begin{definition}
\label{d1}
Let $\kappa \in \mathcal{S}_{-1}$. The $n$-fold GFI ($n \in \N$) is defined as a composition of $n$ GFIs with the kernel $\kappa$:
\begin{equation}
\label{GFIn}
(\I_{(\kappa)}^n\, f)(t) := (\underbrace{\I_{(\kappa)} \ldots \I_{(\kappa)}}_n\, f)(x) = (\kappa^n\, *\, f)(t),\ \kappa^n:= \underbrace{\kappa* \ldots * \kappa}_n,\ t>0.
\end{equation}
\end{definition}

It is worth mentioning that the kernel $\kappa^n,\ n\in \N$ is from the space 
$C_{-1}(0,+\infty)$, but it is not always a Sonine kernel. For example,         in the case of the Riemann--Liouville fractional integral $I_{0+}^\alpha,\ 0<\alpha<1$, the operator \eqref{GFIn} is the Riemann--Liouville fractional integral of the order $n\alpha,\ n\in \N$. Its kernel is $h_{n\alpha}$ that is a Sonine kernel only in the case 
$\alpha<1/n$. However, the operator $I_{0+}^\alpha$ is considered to be a fractional integral for any $\alpha >0$, i.e., also in the case when its kernel $h_{\alpha}$ is not a Sonine kernel and not singular at the point zero. The common justification for this actual situation is that one can define an appropriate fractional derivative, the Riemann--Liouville fractional derivative, that is connected to the Riemann--Liouville fractional integral through the first and second fundamental theorems of FC. In the rest of this section, we introduce an analogous construction for the GFI \eqref{GFI}. 

\begin{definition}
\label{d2}
Let $\kappa \in \mathcal{S}_{-1}$ and $k$ be its associate Sonine kernel. The $n$-fold GFD ($n \in \N$) in the Riemann--Liouville sense is defined as follows:
\begin{equation}
\label{GFDLn}
(\D_{(k)}^n\, f)(t) := \frac{d^n}{dt^n} ( k^n * f)(t),\, k^n:= \underbrace{k* \ldots * k}_n,\ t>0.
\end{equation}
\end{definition}

The $n$-fold GFD is a generalization of the Riemann--Liouville fractional derivative of an arbitrary order $\alpha > 0$ to the case of the general Sonine kernels. For example,         for the Sonine kernels pair $\kappa(t) = h_{\alpha}(t),\ k(t) = h_{1-\alpha}(t),\ 0<\alpha < 1$, the two-fold GFD \eqref{GFDLn} has the kernel $k(t) = h_{2-2\alpha}(t)$ and thus it can be represented in terms of the Riemann--Liouville fractional derivative:
\begin{equation}
\label{RLn}
(\D_{(k)}^2\, f)(t) = (D_{0+}^{2\alpha}\, f)(t) = 
\begin{cases} 
\frac{d^2}{dt^2}(I^{2-2\alpha}_{0+}\, f)(t),& \frac{1}{2} < \alpha <1,\ t>0, \\
 \frac{d}{dt}(I^{1-2\alpha}_{0+}\, f)(t),& 0 < \alpha \le \frac{1}{2},\ t>0. 
 \end{cases}
\end{equation}

Definition \ref{d2} of the the $n$-fold GFD immediately implicates the following recurrent formula: 
\begin{equation}
\label{GFDLn_rec}
(\D_{(k)}^n\, f) (t) = 
\frac{d}{dt}\frac{d^{n-1}}{dt^{n-1}} \left( (k^{n-1}*(k * f))(t)\right) = 
\frac{d}{dt} (\D_{(k)}^{n-1}\, (k*f)) (t),\ n\in \N.
\end{equation}

This formula connects the $n$-fold GFD with the $(n-1)$-fold GFD and is useful for extension of some results derived for the GFD \eqref{GFDL} to the case of the $n$-fold GFD \eqref{GFDLn}. As an example, we generalise the representation \eqref{GFDL-1} of the GFD \eqref{GFDL} to the case of the $n$-fold GFD \eqref{GFDLn}. 

\begin{theorem}
\label{t5}
Let $\kappa \in \mathcal{S}_{-1}$ and $k$ be its associate Sonine kernel. 

Then,  the representation 
\begin{equation}
\label{GFDLn-1}
(\D_{(k)}^n\, f) (t) = (k^n\, * \, f^{(n)})(t) + \sum_{j=0}^{n-1}f^{(j)}(0) 
\frac{d^{n-j-1}}{dt^{n-j-1}} k^n(t),\ t>0
\end{equation}
holds true for any function $f\in C_{-1}^n(0,+\infty)$. 
\end{theorem}

\begin{proof}
We start with a proof of  Formula \eqref{GFDLn-1} in the case $n=2$. The inclusion 
 \mbox{$f\in C_{-1}^2(0,+\infty)$} implies that $f\in C^1[0,+\infty)$ (Theorem \ref{t2}) and $f^{\prime\prime} \in C_{-1}(0,+\infty)$. For $k\in C_{-1}(0,+\infty)$, we immediately get the relations (see the proof of Theorem \ref{t3} for derivation of a similar formula) 
$$
(k * f)(0) = \lim_{t\to 0+} (k * f)(t) = 0,\ (k * f^\prime)(0) = \lim_{t\to 0+} (k * f^\prime)(t) = 0.
$$

Then,  we employ the recurrent Formula \eqref{GFDLn_rec} and the representation \eqref{GFDL-1} for the GFD \eqref{GFDL} and arrive at the following chain of equations:
$$
(\D_{(k)}^2\, f) (t) = \frac{d}{dt} (\D_{(k)}\, (k*f)) (t) 
 = \frac{d}{dt} ((k*(k*f)^\prime)(t) + (k*f)(0)k(t)) =
$$
$$
\frac{d}{dt} ((k*(k*f)^\prime)(t)) = \frac{d}{dt} (k*\D_{(k)}\, f)(t) = \frac{d}{dt} 
(k*((k*f^\prime)(t) + f(0)k(t))) = 
$$
$$
\frac{d}{dt} (k*(k*f^\prime))(t) + \frac{d}{dt} f(0)(k*k)(t) = (\D_{(k)}\, k*f^\prime)(t) + f(0)\frac{d}{dt} k^2(t) = 
$$
$$
(k*(k*f^\prime)^\prime)(t) + (k*f^\prime)(0)k(t)+ f(0)\frac{d}{dt} k^2(t) =
(k*(\D_{(k)}f^{\prime}))(t) + f(0)\frac{d}{dt} k^2(t) =
$$
$$
(k*((k*f^{\prime \prime})(t) + f^{\prime}(0)k(t)))(t) + f(0)\frac{d}{dt} k^2(t) =
$$
$$
(k^2*f^{\prime \prime})(t) + f(0)\frac{d}{dt} k^2(t) + f^{\prime}(0)k^2(t).
$$

The proof of    Formula \eqref{GFDLn-1} for $n=3,4,\dots$ is curried out by employing the method of the mathematical induction and using the recurrent Formula \eqref{GFDLn_rec} and the same technique we employed for the proof in the case $n=2$ based on the formula for $n=1$. 
\end{proof}

As in the case of the GFD \eqref{GFDL}, the representation 
\eqref{GFDLn-1} is used to define a Caputo type $n$-fold GFD. 

\begin{definition}
\label{d3}
Let $\kappa \in \mathcal{S}_{-1}$ and $k$ be its associate Sonine kernel. The $n$-fold GFD ($n \in \N$) in the Caputo sense is defined as follows:
\begin{equation}
\label{GFDCn}
( _*\D_{(k)}^n\, f)(t) := (\D_{(k)}^n\, f) (t) - \sum_{j=0}^{n-1}f^{(j)}(0) 
\frac{d^{n-j-1}}{dt^{n-j-1}} k^n(t),\ t>0.
\end{equation}
\end{definition}

Simple calculations show that the $n$-fold Caputo type GFD \eqref{GFDCn} can be also represented as a regularised $n$-fold GFD in the Riemann--Liouville sense:
\begin{equation}
\label{GFDCn-new}
( _*\D_{(k)}^n\, f)(t) = \left(\D_{(k)}^n\, \left[f(\cdot) - \sum_{j=0}^{n-1}f^{(j)}(0)\, h_{j+1}(\cdot) \right]\right)(t),\ t>0.
\end{equation}

According to Theorem \ref{t5}, on the space $C_{-1}^n(0,+\infty)$, the $n$-fold GFD $ _*\D_{(k)}^n$ takes the following simple form:
\begin{equation}
\label{GFDCn_1}
( _*\D_{(k)}^n\, f)(t) = (k^n\, * \, f^{(n)})(t),\ t>0.
\end{equation}

As in the case of the Caputo fractional derivative \eqref{CD}, the $n$-fold GFD in the form \eqref{GFDCn} makes sense also for the functions that are not $n$-times differentiable (e.g., for the functions from the space $AC^{n-1}[0,+\infty)$). However, the representation of type \eqref{GFDCn_1} for the Caputo fractional derivative is very convenient and thus widely used. 

The denotation ``$n$-fold GFD'' is justified by the fact that the $n$-fold GFDs both in the Riemann--Liouville and in the Caputo senses are left inverse operators to the $n$-fold GFI \eqref{GFIn}. 

\begin{theorem}[First Fundamental Theorem for the $n$-fold GFD]
\label{t6}
Let $\kappa \in \mathcal{S}_{-1}$ and $k$ be its associate Sonine kernel. 

Then,  the $n$-fold GFD \eqref{GFDLn} is a left inverse operator to the $n$-fold GFI \eqref{GFIn} on the space $C_{-1}(0,+\infty)$: 
\begin{equation}
\label{FTLn}
(\D_{(k)}^n\, \I_{(\kappa)}^n\, f) (t) = f(t),\ f\in C_{-1}(0,+\infty),\ t>0
\end{equation}
and the GFD \eqref{GFDCn} is a left inverse operator to the $n$-fold GFI \eqref{GFIn} on the space $C_{-1,k}^n(0,+\infty) := \{f:\ f(t)=(\I_{(k)}^n\, \phi)(t),\ \phi\in C_{-1}(0,+\infty)\}$: 
\begin{equation}
\label{FTCn}
( _*\D_{(k)}^n\, \I_{(\kappa)}^n\, f) (t) = f(t),\ f\in C_{-1,k}^n(0,+\infty),\ t>0.
\end{equation}
\end{theorem}
\begin{proof}
The proof of    Formula \eqref{FTLn} repeats the arguments employed for the proof of Theorem \ref{t3} (the case $n=1$):
$$
(\D_{(k)}^n \I_{(\kappa)}^n\, f) (t) = \frac{d^n}{dt^n}(k^n * (\I_{(\kappa)}^n\, f))(t) = 
\frac{d^n}{dt^n}(k^n * \kappa^n * f)(t) = 
$$
$$
\frac{d^n}{dt^n}(\{ 1\}^n * f)(t) = \frac{d^n}{dt^n} (I_{0+}^n f)(t) = f(t).
$$

For a function $f\in C_{-1,k}^n(0,+\infty)$, we show that
\begin{equation}
\label{cond}
\I_{(\kappa)}^n\, f \in C_{-1}^n(0,+\infty)\ \mbox{and} \ \frac{d^j}{dt^j} \I_{(\kappa)}^n\, f\Bigl|_{t=0} = 0,\ j=0,1,\dots,n-1. \Bigr.
\end{equation}

Indeed, these properties immediately follow from the representation $f(t)=(\I_{(k)}^n\, \phi)(t),$ $\phi\in C_{-1}(0,+\infty)$:
$$
\I_{(\kappa)}^n\, f = \I_{(\kappa)}^n \I_{(k)}^n\, \phi = \kappa^n*k^n*\phi = \{ 1\}^n*\phi, \ \phi\in C_{-1}(0,+\infty).
$$

Then,  we employ    Formula \eqref{FTLn} and arrive at the final result: 
$$
( _*\D_{(k)}^n\, \I_{(\kappa)}^n\, f) (t) = (\D_{(k)}^n\, \I_{(\kappa)}^n\, f) (t) - \sum_{j=0}^{n-1}\frac{d^j}{dt^j} \I_{(\kappa)}^n\, f\Bigl|_{t=0}
\frac{d^{n-j-1}}{dt^{n-j-1}} k^n(t) = \Bigr.
$$
$$
(\D_{(k)}^n\, \I_{(\kappa)}^n\, f) (t) = f(t).
$$
\end{proof}

\begin{remark}
\label{r2}
The same arguments and derivations that we employ  in Remark \ref{r1} for the space $C_{-1,k}(0,+\infty)$ are applicable for the space 
$$
C_{-1,k}^n(0,+\infty) = \{f:\ f(t)=(\I_{(k)}^n\, \phi)(t),\ \phi\in C_{-1}(0,+\infty)\}.
$$ 

As a result, this space can be also characterised as follows:
$$
C_{-1,k}^n(0,+\infty) = \{f:\ \I_{(\kappa)}^n f \in C_{-1}^n(0,+\infty)\, \wedge \, \frac{d^j}{dt^j} \I_{(\kappa)}^n\, f\Bigl|_{t=0} = 0,\ j=0,1,\dots,n-1\}. \Bigr.
$$
\end{remark}
 
Finally, we formulate and prove the second  fundamental theorem for the $n$-fold GFDs in the Caputo sense. 

\begin{theorem}[Second Fundamental Theorem for the $n$-fold GFD]
\label{t7}
Let $\kappa \in \mathcal{S}_{-1}$ and $k$ be its associate Sonine kernel. 

Then,  the relation
\begin{equation}
\label{2FTCn}
(\I_{(\kappa)}^n\, _*\D_{(k)}^n\, f) (t) = f(t)-\sum_{j=0}^{n-1} \, f^{(j)}(0)\frac{t^{j}}{j!},\ t>0
\end{equation}
holds true for any function $f\in C_{-1}^n(0,+\infty)$.
\end{theorem}

\begin{proof}
For the functions 
$f\in C_{-1}^n(0,+\infty)$, the GFD $ _*\D_{(k)}^n$ can be represented in the form \eqref{GFDCn_1}. Moreover, according to Theorem \ref{t2}, the inclusion $f\in C_{-1}^n(0,+\infty)$ implies the inclusion $f\in C^{n-1}[0,\, +\infty)$. Then,  we have the following chain of equations:
$$
(\I_{(\kappa)}^n\, _*\D_{(k)}^n\, f) (t) = (\kappa^n * ( _*\D_{(k)}^n\, f))(t) = (\kappa^n * k^n*f^{(n)})(t) = 
$$
$$
(\{ 1\}^n *f^{(n)})(t) = f(t)-\sum_{j=0}^{n-1} \, f^{(j)}(0)\frac{t^{j}}{j!},\ t>0.
$$
\end{proof}
For the $n$-fold GFD in the Riemann--Liouville sense, under appropriate conditions on its kernel $k$, we get the formula
\begin{equation}
\label{2FTLn}
(\I_{(\kappa)}^n\, \D_{(k)}^n\, f) (t) = f(t),\ f\in C_{-1}^n(0,+\infty),\ t>0.
\end{equation}
In particular, if the conditions
\begin{equation}
\label{2FTLn-c}
(\kappa^n(\tau) * \frac{d^{n-j-1}}{d\tau^{n-j-1}} k^n(\tau))(t) = 
\frac{d^{n-j-1}}{dt^{n-j-1}} (\kappa^n(\tau) * k^n(\tau))(t),\ j=0,\dots,n-1
\end{equation}
are satisfied, then    Formula \eqref{2FTLn} holds true:
$$
(\I_{(\kappa)}^n\, \D_{(k)}^n\, f) (t) = (\I_{(\kappa)}^n\, ((k^n*f^{(n)})(\tau) + \sum_{j=0}^{n-1}f^{(j)}(0) 
\frac{d^{n-j-1}}{d\tau^{n-j-1}} k^n(\tau))) (t) =
$$
$$
(\I_{(\kappa)}^n\, _*\D_{(k)}^n\, f) (t) + (\kappa^n * (\sum_{j=0}^{n-1}f^{(j)}(0) 
\frac{d^{n-j-1}}{d\tau^{n-j-1}} k^n(\tau)))(t) = 
$$
$$
f(t)-\sum_{j=0}^{n-1} \, f^{(j)}(0)\frac{t^{j}}{j!} + \sum_{j=0}^{n-1} \, f^{(j)}(0) 
\frac{d^{n-j-1}}{dt^{n-j-1}} (\kappa^n(\tau) * k^n(\tau))(t) = 
$$
$$
f(t)-\sum_{j=0}^{n-1} \, f^{(j)}(0)\frac{t^{j}}{j!} + \sum_{j=0}^{n-1} \, f^{(j)}(0) 
\frac{d^{n-j-1}}{dt^{n-j-1}}(\{ 1\}^n(t)) =
$$
$$
f(t)-\sum_{j=0}^{n-1} \, f^{(j)}(0)\frac{t^{j}}{j!} + \sum_{j=0}^{n-1} \, f^{(j)}(0)\frac{t^{j}}{j!}
= f(t).
$$

Evidently, the conditions \eqref{2FTLn-c} are fulfilled, for instance, in the case $n=1$ (Theorem \ref{t4}) and for the kernel $k(t)=h_{1-\alpha}(t),\ 0<\alpha < 1/n$ of the Riemann--Liouville fractional derivative \eqref{RLD} and then   Formula \eqref{2FTLn} is valid for $f\in C_{-1}^n(0,+\infty)$. 

\section{Conclusions}
\label{sec:5}
Starting from the publication \cite{Koch11}, a considerable interest in FC  has been devoted to the fractional integrals and derivatives with the general kernels. Because the GFDs are the left-inverse operators to the corresponding GFIs (they have to satisfy the first fundamental theorem of FC), their kernels cannot be arbitrary. As shown in \cite{DGGS}, a necessary condition that these kernels have to fulfil is their singularity at the point zero. In \cite{Han20}, the Sonine condition was introduced as a sufficient condition for the kernels of the GFIs and GFDs that satisfy the first fundamental theorem of FC. Because the Sonine kernels are singular at the point zero, they satisfy the necessary condition formulated in \cite{DGGS}. Thus, the operators with the Sonine kernels are of particular interest in the framework of FC. In \cite{Koch11,Han20}, the GFIs and GFDs with the Sonine kernels from the sets $\mathcal{K}$ and $\mathcal{H}$ (see Section \ref{sec:3}), respectively,  are  considered. 

In this paper, we introduce  and study  the GFIs and GFDs with the Sonine kernels from an essentially larger set compared to $\mathcal{K}$ and $\mathcal{H}$, namely  from the set $\mathcal{S}_{-1}$ (see \mbox{Section \ref{sec:3}} for its definition). The kernels from $\mathcal{S}_{-1}$ are continuous functions that possess integrable singularities of the power function type at the point zero. Some basic properties of the GFIs and GFDs including the first and second fundamental theorems of FC  are  deduced on the spaces of functions with an integrable singularity at the point zero. These spaces are very natural while studying the GFIs and GFDs with the Sonine kernels from $\mathcal{S}_{-1}$. 

The GFDs treated  so  far are  generalisations of the Riemann--Liouville and the Caputo fractional derivatives of the order $\alpha$ restricted to the interval $(0,\ 1)$. In this paper, we introduce  the $n$-fold GFDs that can be considered as generalisations of the Riemann--Liouville and the Caputo fractional derivatives of an arbitrary positive order $\alpha$. For these operators, the first and second fundamental theorems of FC as well as some other properties  are  formulated and proved. 

In this paper, we focus  on the mathematical properties of the GFIs and GFDs with the Sonine kernels from $\mathcal{S}_{-1}$. However, these operators are also potentially useful for applications. One of the main problems in applications of the FC operators and the fractional differential equations is to decide which kind of the fractional derivatives should be chosen by ``fractionalisation'' of the conventional models containing the differential operators. The models with the GFDs with the Sonine kernels from $\mathcal{S}_{-1}$ provide a wider scope for fitting the data available from the measurements compared to the models with the conventional fractional derivatives. Moreover, up to now, the fractional differential equations with a single GFD have been considered in the literature. The $n$-fold GFDs introduced in this paper open a way for formulation of the models in form of the multi-term fractional differential equations with several GFDs. These equations are a far reaching generalisation of the multi-term fractional differential equations with the Riemann--Liouville or Caputo fractional derivatives that are actively employed for modelling of the multi-scale processes. These equations and their potential applications will be considered elsewhere.

\end{document}